\documentstyle{amsppt}
\magnification=1200

\hoffset=-0.5pc
\vsize=57.2truepc
\hsize=38truepc
\nologo
\spaceskip=.5em plus.25em minus.20em

\define\Bobb{\Bbb} 
\define\fra{\frak} 
 \define\abramars{1}
 \define\armcusgo{2}
 \define\armgojen{3}
 \define\arnobook{4} 
 \define\atibottw{5}
 \define\berghueb{6}
 \define\bocosroy{7}
\define\brymcone{8}
\define\goldmone{9}                                                            
\define\goldmtwo{10}
\define\golmilon{11}
\define\golmiltw{12}
\define\gormacon{13}
\define\guhujewe{14}
\define\hartsboo{15}
 \define\howeone{16}
\define\poiscoho{17}
 \define\souriau{18}
 \define\singula{19}
\define\singulat{20}
  \define\smooth{21}
 \define\poisson{22}
 \define\locpois{23}
   \define\modus{24}
\define\modustwo{25}
  \define\moscou{26}
   \define\srni{27}
\define\claustha{28}
\define\kan{29}
\define\poismod{30}
\define\kaehler{31}
\define\huebjeff{32}
\define\kapmione{33}
\define\kemneone{34}
\define\kobanomi{35}
\define\kunzbotw{36}
\define\lermonsj{37}
\define\marswein{38}
\define\naramntw{39}
\define\naramnth{40}
\define\narashed{41}
\define\newstboo{42}
\define\gwschwar{43}
\define\gwschwat{44}
\define\seshaone{45}
\define\seshaboo{46}
\define\sjamlerm{47}
 \define\weiltwo{48}
\define\weinstwo{49}
\define\weinsone{50}
\define\weinsthi{51}
\define\weylbook{52}
\define\whitnone{53}
\define\whitnboo{54}
\topmatter
\title 
Singularities and Poisson geometry of certain representation spaces
\endtitle
\author Johannes Huebschmann
\endauthor
\affil 
Universit\'e des Sciences et Technologies de Lille
\\
UFR de Math\'ematiques
\\
59 655 VILLENEUVE D'ASCQ C\'edex
\\
Johannes.Huebschmann\@agat.univ-lille1.fr
\endaffil
\date{September 26, 2000}\enddate
\endkeywords
\endsubjclass
\endtopmatter
\document
\leftheadtext{Johannes Huebschmann}

\beginsection 1. Introduction

Certain representation spaces have been investigated
by algebraic geometers as moduli spaces of
holomorphic bundles over a Riemann surface.
Such moduli spaces exhibit symplectic
and K\"ahler structures 
as well as gauge theory interpretations.
The purpose of this article is
to elucidate the local structure
of such a space,
and the focus will be on the singularities.
Among the tools will be the interconnection between
the theory of algebraic 
and symplectic
quotients
and, furthermore,
Poisson structures,
a concept which
has been  known in mathematical physics 
for long and is
currently of much interest
in mathematics as well.

On 
$\Bobb R^{2n}$ with its standard coordinates
$q_1,\dots,q_n,p_1,\dots,p_n$, the formula
$$
\{f,h\} = \sum\left(\frac{\partial f}{\partial p_j}
\frac{\partial  h}{\partial q_j}
-
\frac{\partial f}{\partial q_j}
\frac{\partial  h}{\partial p_j}
\right)
\tag 1.1
$$
yields
a bracket $\{\cdot,\cdot\}$
on the algebra of smooth functions,
nowadays referred to as a Poisson bracket.
This
bracket was introduced by {\smc Poisson\/}
around 1809, and he observed that,
given three functions $f,g,h$ with
$
\{f,g\}= 0
$
and
$
\{f,h\}= 0,
$
one also has
$
\{f,\{g,h\}\}= 0.
$
This means that if
$g$ and $h$ are integrals of motion
for
(the hamiltonian vector field of) $f$,
so is 
$\{g,h\}$.
See e.~g. \cite{\abramars\  p.~196},
\cite{\arnobook\  p.~216}.
Abstractly, 
a Poisson algebra is a commutative algebra $A$ together
with the additional structure
of a Lie bracket $\{\ ,\ \}$ which behaves as a derivation in each variable
with respect to the algebra structure.
A smooth Poisson structure on an ordinary 
smooth manifold
is symplectic if and only if
it is locally of the kind (1.1).

Poisson structures provide
some insight into the geometry
of certain representation spaces.
Thus, let $G$ be a compact Lie group and  $\pi$ 
the fundamental group of a closed surface $\Sigma$, 
and consider the {\it representation space\/}
$\roman{Rep}(\pi,G) = \roman{Hom}(\pi,G)\big / G$;
here $G$ acts on $\roman{Hom}(\pi,G)$
via conjugation.
Pick an invariant symmetric non-degenerate bilinear form on the Lie
algebra $\fra g$ of $G$.
The following result
may be found in \cite\modus; see also
\cite{\moscou--\claustha} for a leisurely introduction
and for more references. 

\proclaim{Theorem 1} 
The decomposition of
$\roman{Rep}(\pi,G)$
into orbit types of representations is a stratification,
and the data induce a
Poisson algebra 
$$
\left(C^{\infty}(\roman{Rep}(\pi,G)),\{\ ,\ \}\right)
$$
of continuous functions on $\roman{Rep}(\pi,G)$
turning
the latter
into a stratified symplectic space in such a way that
a choice of complex structure on $\Sigma$
induces a K\"ahler structure on each stratum.
Furthermore, the space
$\roman{Rep}(\pi,G)$ is locally semi-algebraic.
Finally, the Poisson structure detects the stratification.
\endproclaim

Here 
the algebra
$C^{\infty}(\roman{Rep}(\pi,G))$
encapsulates 
the {\it real\/} geometry of the space 
$\roman{Rep}(\pi,G)$.
With this structure,
$\roman{Rep}(\pi,G)$ is {\it not\/} a real variety in the usual sense;
its singularity behaviour
is alluded to
by the wording \lq\lq locally semi-algebraic\rq\rq.
The meaning of \lq\lq (locally)
semi-algebraic\rq\rq\ will be explained in Section 2 below
where we will also explain how the
Poisson structure detects the stratification,
and the decomposition into orbit types will be reproduced in Section 4
below.

There is a version of the above result
involving twisted  representation spaces
\cite\modus\ 
as well as one involving surfaces with non-empty boundary
\cite\guhujewe.
In algebraic geometry,
after a choice of complex structure
on $\Sigma$ has been made,
for $G=\roman U(n)$, the space
$\roman{Rep}(\pi,G)$
arises as the 
(real space underlying the)
moduli space
of semi-stable holomorphic
rank $n$ vector bundles on $\Sigma$ of degree zero,
and 
the twisted  representation spaces
correspond to the case of arbitrary degree
\cite{\naramntw--\newstboo,\, \seshaone,\,
\seshaboo}.
The Poisson structure then gives rise to some interesting Poisson geometry.
In particular,
for
$G=\roman{SU}(2)$,
the space
$\roman{Rep}(\pi,G)$
amounts to the moduli space of semi-stable holomorphic
vector bundles on $\Sigma$
of rank 2, degree 0, and trivial determinant. 
To describe the stratification in this case, let
$Z \subseteq G$ denote the center
and $T \subseteq G$
a maximal torus.
The space $N = \roman{Rep}(\pi,G)$ is a union  
$$
N = N_G \cup N_{(T)}\cup N_Z
$$
of three strata
where $N_{(K)}$
denotes the points of orbit type $(K)$.
The piece  $N_Z$
is called the top stratum.

For genus $\ell \geq 2$,
{\smc Narasimhan-Ramanan}
proved that the complement
$\Cal K$ of the top stratum
is the {\it Kummer\/} variety
of $\Sigma$
associated with its Jacobian $\Cal J$
and the canonical involution thereupon
\cite\naramntw.
This has the following consequence,
established in
\cite\locpois.

\proclaim{Theorem 2}
For $G=\roman{SU}(2)$,
when $\Sigma$ has genus $\ell \geq 2$,
the Poisson algebra
$(C^{\infty}(\roman{Rep}(\pi,G)),\{\cdot,\cdot\})$
detects the Kummer variety
$\Cal K$ in
$\roman{Rep}(\pi,G)$
together with its $2^{2\ell}$ double points.
More precisely,
$\Cal K$ consists of the points
where the rank of the Poisson structure
is not maximal,
the double points being those
where the rank is zero.
\endproclaim

For genus $\geq 3$,
the Kummer variety
$\Cal K$
is precisely the
(complex analytic)
singular locus of $\roman{Rep}(\pi,G)$,
a result due to
{\smc Narasimhan-Ramanan}
\cite\naramntw.
This has been reproved in \cite\locpois\ 
within our framework.
When $\Sigma$ has genus two the space
$\roman{Rep}(\pi,G)$
equals complex projective 3-space
\cite\naramntw\ 
and
$\Cal K$
is the Kummer surface associated with the Jacobian
of $\Sigma$.
In the literature, this case has been considered
somewhat special:
The algebra of ordinary smooth functions 
on 
$\roman{Rep}(\pi,G)$, realized as
complex projective 3-space,
is a smooth structure 
on $\roman{Rep}(\pi,G)$
in the above sense as well---call it the {\it standard\/} structure---and
the Kummer surface is certainly {\it not\/} the singular locus for this
structure.
However, from our point of view there is {\it no\/}
exception.
Our smooth structure 
$C^{\infty}(\roman{Rep}(\pi,G))$
is {\it not\/}
the standard one,
and 
even in the genus two case,
as a 
stratified symplectic space,
$\roman{Rep}(\pi,G)$
still has singularities:
The Poisson algebra
$(C^{\infty}(\roman{Rep}(\pi,G)),\{\cdot,\cdot\})$
detects 
(the real space underlying)
a Kummer
surface together with its 16 singularities and hence
the underlying algebra
of functions can plainly {\it not\/} be that
of ordinary smooth functions;
in particular, the symplectic structure
on the top stratum
does {\it not\/}
extend to the whole space.
It is interesting to observe that the stratification
mentioned in Theorem 1 is {\it finer\/}
than the standard complex analytic one on complex
projective 3-space.

Under the circumstances of Theorem 1,
the K\"ahler structures on the strata of $\roman{Rep}(\pi,G)$
fit together in a very precise way.
Abstracting from this situation, we isolated the notion of
{\it stratified K\"ahler space\/} in \cite\kaehler.
The structure on complex
projective 3-space
just explained is a particular case of
a stratified K\"ahler structure
which is {\it not\/} the standard K\"ahler structure.
In \cite\kaehler\ it is shown that, in particular,
such \lq\lq exotic\rq\rq\ 
structures on
complex projective spaces
abound.

I am indebted to 
W. Goldman, J. Stasheff, and
the referee
for a number of remarks
which helped improve the
exposition.

\medskip\noindent{\bf 2. Singularities, real vs complex; an (almost toy) 
example} \smallskip\noindent
What is a singularity?
\lq\lq A singular point is one which behaves differently from
the other points close to it.\rq\rq\ 

For illustration
consider 
$\Bobb R^2$, but viewed as a half-cone.
What kind of structure can we put on
$\Bobb R^2$
to distinguish the cone point, $o$, from the other points?
A possible answer: Consider the ordinary coordinates $x,y$ in the plane
and introduce a third variable $\rho$, subject to the relation
$$
x^2 + y^2 = \rho^2.
\tag2.1
$$
Let $A$ be the algebra of smooth functions in these variables,
subject to the relation (2.1).
Notice that $\rho$ is {\it not\/} a smooth function on the plane in the
usual sense, and the algebra $A$ is strictly larger than that
of ordinary smooth functions on the plane.
The situation is entirely parallel to that of an algebra of continuous 
functions on complex projective 3-space
strictly larger than the ordinary one
mentioned in Theorem 2 above.
Under the present circumstances,
the real algebraic variety $Z$ defined by (2.1),
that is, the set of points $(x,y,\rho)$ satisfying this equation,
is a (double) cone. The half cone we are interested in is characterized
by the additional constraint $\rho \geq 0$.
Technically one says that
the half cone $C$ is a {\it semi-algebraic\/} subset of $Z$,
cf. e.~g. \cite\bocosroy.
Why do we distinguish $C$ from $Z$?
The answer is that a space of the kind
$\roman{Rep}(\pi,G)$
is locally semi-algebraic.

Given a space $X$ which is decomposed into pieces which
are smooth manifolds,
a {\it smooth\/} structure on $X$
is an algebra $C^{\infty}(X)$ of 
(real or complex valued)
continuous functions which, on each piece,
restrict to ordinary smooth functions.
A space $X$ together with a smooth structure
is called a {\it smooth\/} space.

The half cone $C$, together with the algebra $A=C^{\infty}(C)$ 
of functions defined above, is a smooth space, and so is 
the (double) cone $Z$, together with $A=C^{\infty}(Z)$.
Viewed abstractly,
$Z$ is just a realization of the real maximal spectrum
of $A$ and,
as algebras,
$C^{\infty}(C)$
and
$C^{\infty}(Z)$
coincide.
We can view 
$A=C^{\infty}(C)=C^{\infty}(Z)$
as the algebra of ordinary smooth functions on
$\Bobb R^3$, subject to the relation (2.1),
which means that every function belonging to this smooth structure
is the restriction of an ordinary smooth function
on the ambient space. 
More generally: Given a closed connected subset
$Y \subseteq \Bobb R^n$ (not necessarily a smooth manifold),
the constinuous functions 
on $Y$ which are restrictions of ordinary smooth functions
on $\Bobb R^n$
endow $Y$ with a smooth structure;
this construction goes back to {\smc Whitney\/}
\cite\whitnone, and
the resulting functions on $Y$ are nowadays
called {\it Whitney\/}-smooth functions.
Thus, in the present situation,
the algebra
$A=C^{\infty}(C)=C^{\infty}(Z)$
is that of Whitney smooth functions on
$Z$ as well as on $C$.

Given a smooth space
$(X,C^{\infty}(X))$,
for each point $x$ of $X$, the {\it ideal\/}
$\fra m_x$ of $x$ in $C^{\infty}(X)$
consists of the smooth functions which vanish at $x$.
See e.~g. \cite{\gwschwar,\, \gwschwat}.
The {\it Zariski tangent space\/}
$\roman T_xX$ of $X$ at $x$ is the vector space
$(\fra m_x\big /\fra m^2_x)^*$ \cite{\bocosroy,\, \whitnboo};
see e.~g. \cite\hartsboo\ for the notion of Zariski tangent space
in algebraic geometry.

Another well known description
of the 
Zariski tangent space is this:
Let $x \in X$ and
view $\bold R$ as a
$C^{\infty}(X)$-module,
written
$\bold R_x$,
by means of
the evaluation mapping
from $C^{\infty}(X)$ to $\bold R$
which assigns to a function
$f$ its value $f(x)$ at $x \in X$;
now
a {\it derivation\/}
at $x \in X$
is a linear map
$d$ from $C^{\infty}(X)$ to $\bold R$
satisfying the usual {\it Leibniz\/} rule
$$
d(fh) = (df)h(x) + f(x)dh .
$$
We denote the real vector space
of all 
derivations
at $x$
by 
$\roman{Der}(C^{\infty}(X),\bold R_x)$.
For $x \in X$, the assignment
to $\phi \in \roman T_x X$
of the derivation
$d_{\phi}$ at $x$
given by
$d_{\phi} (f) = \phi( f - f_x)$
identifies 
$\roman T_x X$
with 
$\roman{Der}(C^{\infty}(X),\bold R_x)$;
here $f \in C^{\infty}(X)$ and
$f_x$ denotes the  function 
having constant value
$f(x)$.
Thus, when $X$ is a smooth manifold near a point $x$
in the usual sense,
with standard smooth structure near $x$,
the Zariski tangent space boils down to the ordinary smooth
tangent space
whence there is no risk of confusion in notation.
Given smooth spaces $(X,C^{\infty}(X))$,
$(Y,C^{\infty}(Y))$,
and a smooth  map $\phi$ from $X$  to $Y$,
the {\it derivative\/}
at a point $x \in X$
is the dual
${
d\phi_x
\colon
\roman T_x X
@>>>
\roman T_{\phi x} Y
}$
of the linear map
from
$\bold m_{\phi(x)}\big/\bold m^2_{\phi(x)}$ 
to
$\bold m_x\big/\bold m^2_x$
induced by $\phi$.

Intuitively,
the Zariski tangent space at $x$ is the linear span
of all the {\it tangents\/} to $x$.
Lack of space does not allow us to elaborate on this here;
see e.~g. \cite\kunzbotw.
In the above example of a (half) cone:
At every point different from $o$, the Zariski tangent space
is a plane while, at $o$, it is a copy of $\Bobb R^3$.
In fact, this copy of
$\Bobb R^3$
is the span of the cone.
Thus we may distinguish the singular point $o$
from the other points by means of the notion of Zariski tangent space
defined in terms of the
{\it non-standard smooth structure\/}
involving the additional function $\rho$.

There is yet another way to distinguish the cone point from the other 
points: Introduce a Poisson structure $\{\ ,\ \}$ by setting
$$
\{x,y\} = 2\rho,
\quad
\{x,\rho\} = 2y,
\quad
\{y,\rho\} = -2x.
\tag2.2
$$
This Poisson structure is symplectic
everywhere except at the cone point
where it degenerates.
Thus it distinguishes the cone point from the other points.

Recall that a {\it symplectic\/} manifold is a smooth manifold $M$ together
with a closed non-degenerate 2-form $\omega$.
Given $f$, the identity $\omega(X_f,\cdot ) = df$ then associates a uniquely
determined vector field $X_f$ to $f$, the {\it Hamiltonian vector field\/}
of $f$, and the assignment $\{f,h\} = X_f h$
yields a Poisson bracket on the  algebra $C^{\infty} (M)$
of ordinary smooth functions on $M$.
Such a Poisson bracket is called {\it symplectic\/}.
The Poisson structure (1.1) is a special case thereof.

A decomposition of a space $X$ into 
pieces which are
smooth manifolds
such that these pieces fit together in a certain precise way
is called a {\it stratification\/}
\cite\gormacon.
More precisely:
Let $X$ be a Hausdorff paracompact topological space
and let $\Cal I$ be a partially ordered set
with order relation denoted by $\leq$.
An $\Cal I$-decomposition of $X$ is a locally finite collection
of disjoint locally closed manifolds $S_i \subseteq X$
called {\it pieces\/} 
(recall that a collection $\Cal A$ of subsets of $X$ is said to be
{\it locally finite\/} provided every $x \in X$ has a neighborhood
$U_x$ in $X$ such that $U_x \cap A \ne \emptyset$
for at most finitely many $A$ in $\Cal A$)
such that the following hold:
\smallskip
$X = \cup S_i$
\smallskip
$S_i\cap \overline S_j \ne \emptyset \Longleftrightarrow S_i \subseteq 
\overline S_j \Longleftrightarrow i \leq j$.
\smallskip\noindent
The space $X$ is then called a {\it decomposed\/} space.
A decomposed space $X$ is called a {\it stratified space\/}
if the pieces of $X$, called {\it strata\/},
satisfy
the following condition:
Given a point $x$ in a piece $S$ there is an open neighborhood $U$
of $x$ in $X$,
an open ball $B$ around $x$ in $S$, a stratified space $L$, called the 
{\it link\/} of $x$, and a decomposition preserving homeomorphism
from $B \times C^{\circ}(L)$ onto $U$.
Here $C^{\circ}(L)$
refers to the open cone on $L$ and, as a stratified space,
$L$ is less complicated than $C^{\circ}(L)$
whence the definition is not circular;
the idea of complication is here made precise by means of the notion
of {\it depth\/}.
A {\it stratified symplectic space\/}
\cite\sjamlerm\ 
is a stratified space $X$
together with a 
Poisson algebra $(C^{\infty}(X),\{\ ,\ \})$
of continuous functions which, on each piece,
restricts to an ordinary smooth
symplectic Poisson structure;
in particular,
$C^{\infty}(X)$ is a smooth structure on $X$.

The half cone, endowed with the
above smooth structure and Poisson algebra
(2.2) is an example of a stratified symplectic space
with two strata,
the cone point and its complement, the \lq\lq top\rq\rq\ stratum.
The Poisson structure detects the stratification:
It has \lq\lq rank\rq\rq\  zero at the cone point and 
\lq\lq rank\rq\rq\ two
on the top stratum whence it is symplectic there.
Given a point $x$ of $X$,
the {\it rank\/} at $x$ refers here to the rank of the 
linear map
from $\Omega_x(X)$ to $\roman{Der}(C^{\infty}(X), \Bobb R_x)$
induced by the Poisson structure;
more precisely,
the canonical map from $\Omega(X) \otimes_{C^{\infty}(X)}\Bobb R_x$
to $\Omega_x(X)$ ($\cong \fra m_x\big /\fra m^2_x$, see what
is said above,) 
is an isomorphism of vector spaces,
the Poisson structure
induces a $C^{\infty}(X)$-morphism
from   the 
$C^{\infty}(X)$-module $\Omega(X)$ of differentials
with respect to 
$C^{\infty}(X)$, endowed with the
$(\Bobb R,C^{\infty}(X))$-Lie algebra structure
coming from the Poisson structure, cf. \cite\poiscoho,
to 
the $C^{\infty}(X)$-module $\roman{Der}(C^{\infty}(X))$
of derivations of $C^{\infty}(X)$,
and this morphism, in turn,
induces
the linear map
from $\Omega_x(X)$ to $\roman{Der}(C^{\infty}(X), \Bobb R_x)$
under discussion.
More details about the detection
of the stratification by means of the notion of rank
may be found in \cite\locpois.

The ordinary complex
structure of the plane,
combined with the symplectic structure
on the top stratum, turns the latter into a K\"ahler manifold.
But the cone point is not a singularity for the complex structure
whence the latter cannot be used to detect the cone point.
This reflects the remark made earlier
that the symplectic stratification
is finer than the standard complex analytic one.

For completeness, we recall at this stage that a {\it  K\"ahler manifold\/}
is a smooth complex manifold 
together with a hermitian metric whose imaginary part is a symplectic 
structure. Equivalently,
a smooth symplectic
manifold $(M,\omega)$ 
together with a complex structure $J$ 
which is compatible with $\omega$ (i.~e. $\omega(JX,JY) 
=\omega(X,Y) $)
is called a K\"ahler manifold provided
the (real) symmetric bilinear form $g$ given by
$g(X,Y) = \omega(X, JY)$
is positive definite.
See e.~g. \cite\kobanomi\ for details.

\beginsection 3. Symplectic reduction in a nutshell

Given a compact Lie group $G$,
a {\it hamiltonian\/} $G$-space is a smooth symplectic $G$-manifold 
$(M,\omega)$
together with a smooth $G$-equivariant map $\mu$ from $M$ to the dual
$\fra g^*$ of the Lie algebra $\fra g$ of $G$ satisfying the formula
$$
\omega(X_M,\cdot) = X \circ d \mu
\tag3.1
$$
for every $X \in \fra g$; here $X_M$ denotes the vector field on $M$
induced by $X \in \fra g$ via the $G$-action,
and $X$ is viewed as a linear map on $\fra g^*$.
The map $\mu$ is called a {\it momentum mapping\/}
(or {\it moment map\/}).
We  recall that
(3.1) says that,
given $X \in \fra g$,
the vector field $X_M$ is the hamiltonian vector field for 
the smooth function
$X\circ \mu$ on $M$.
See \cite{\abramars, \arnobook} for details.
Given a hamiltonian $G$-space $(M,\omega,\mu)$,
the space $M_{\roman{red}} = \mu^{-1}(0)\big / G$
is called the {\it reduced space\/}
or {\it symplectic quotient\/}.
It carries the smooth structure
$C^{\infty}(M_{\roman{red}}) = (C^{\infty}(M))^G)\big/I^G$,
where
$I^G$
refers to the ideal in the algebra
$(C^{\infty}(M))^G$
of smooth $G$-invariant functions on $M$
which vanish on the zero locus $\mu^{-1}(0)$.
Arms-Cushman-Gotay \cite\armcusgo\ 
established the fact that the symplectic Poisson algebra
on $C^{\infty}(M)$
descends to a Poisson structure 
$\{\ ,\ \}_{\roman{red}}$
on
$C^{\infty}(M_{\roman{red}})$, and
Sjamaar-Lerman \cite\sjamlerm\ have shown that
the decomposition of
$M_{\roman{red}}$ according to orbit types is a stratification
in such a way that 
$(M_{\roman{red}}, C^{\infty}(M_{\roman{red}}),\{\ ,\ \}_{\roman{red}})$
is a stratified symplectic space.
When $M_{\roman{red}}$ is smooth, i.~e. has a single stratum,
the structure comes down to that of a smooth symplectic manifold
on $M_{\roman{red}}$,
the ordinary Marsden-Weinstein reduced space \cite\marswein.

\medskip\noindent
{\bf 4. Representation spaces}\smallskip\noindent
Denote the genus of the closed surface $\Sigma$ 
under consideration
by $\ell$.
The standard presentation
$$
\langle
x_1,y_1,\dots, x_{\ell}, y_{\ell}; r
\rangle
$$
of $\pi=\pi_1(\Sigma)$ where $r = [x_1,y_1] \dots [x_\ell,y_\ell]$
determines a {\it word map\/}
$r\colon G^{2 \ell} @>>> G$,
and 
the choice of generators 
identifies
$\roman{Hom}(\pi,G)$ 
with the subset $r^{-1}(e)$ of
$G^{2 \ell}$
(where $e \in G$ is the neutral element)
whence
$\roman{Rep}(\pi,G)$ may be realized as 
the quotient 
$$
\roman{Rep}(\pi,G)
=
r^{-1}(e)\big/ G.
$$
The situation is very similar to that of a hamiltonian $G$-space.
Indeed a corresponding hamiltonian $G$-space,
called {\it extended moduli space\/},
has been constructed
in \cite{\modus, \huebjeff}.
The extended moduli space then determines the
requisite smooth structure,
and symplectic reduction 
yields the stratified symplectic space structure
on
$\roman{Rep}(\pi,G)$
coming into play in Theorem 1.

The decomposition which underlies the stratification 
of $\roman{Rep}(\pi,G)$
is that by orbit types.
For the reader's convenience, we explain this briefly:
Write $N= \roman{Rep}(\pi,G)$ and,
for a closed subgroup $K \subseteq G$
let
$N_{(K)} \subseteq N$
denote the subspace of
classes $[\phi]$ of
representations $\phi$ having
stabilizer $Z_{\phi}$
conjugate to $K$.
Thus the representation space
$N$ 
decomposes into a disjoint union
of {\it orbit types\/} $\Cal N_{(K)},\,(K)\in \Cal I$,
the indexing set
$\Cal I$ being that of all possible
stabilizer subgroups
modulo conjugacy,
and each orbit type decomposes further into its connected components.
Sometimes 
the decomposition 
into connected components
is the more appropriate one.

\beginsection 5. The local structure

Consider a homomorphism $\phi \colon \pi \to G$.
It represents a point $[\phi]$ of the space
$\roman{Rep}(\pi,G)$.
By means of $\phi$ we endow $\fra g$ with a
$\pi$-module structure.
The Lie bracket on $\fra g$
induces a graded Lie bracket
$[\ ,\ ]_{\phi}$
on the group cohomology group
$\roman H_\phi^*=\roman H^*(\pi, \fra g_{\phi})$,
and the chosen invariant symmetric bilinear form on $\fra g$,
combined with the isomorphism
from $\roman H^2(\pi, \Bobb R)$ onto $\Bobb R$
given by a choice of fundamental class of $\pi$,
induces a non-degenerate graded commutative
pairing
$$
\roman H_\phi^*
\otimes
\roman H_\phi^{2-*}
@>>>
\Bobb R .
\tag5.1
$$
In degree 1, (5.1) amounts
to a symplectic structure
$\sigma_{\phi}$ on the vector space
$\roman H_\phi^1$.
Moreover, 
$\roman H_\phi^0$ is the Lie algebra of the stabilizer
$Z_{\phi} \subseteq G$
of $\phi$, and the assignment
$\Theta_{\phi}(\eta) = \frac 12[\eta,\eta]_{\phi}$
($\eta \in \roman H_\phi^1$)
yields a momentum mapping 
$\Theta_{\phi}\colon\roman H_\phi^1 \to \roman H_\phi^2$ 
for the action 
of $Z_{\phi}$
on $\roman H_\phi^1$
which is therefore
hamiltonian;
here
$\roman H_\phi^2$
is identified with the dual of
$\roman H_\phi^0$ by means of (5.1).
We then have the following result,
given in Theorem 6.3 of \cite\smooth;
for the analogous gauge theory situation,
it has been spelled out in \cite\singula. 

\proclaim{Theorem 3} 
The reduced space
$(\roman H_\phi^1)_{\roman{red}}$,
with its stratified symplectic space structure,
is a local model for 
$\roman{Rep}(\pi,G)$
near $[\phi]$
as a stratified symplectic space.
More precisely, the choice of $\phi$
(in its class $[\phi]$)
induces a diffeomorphism
of an open neighborhood
$W_\phi$ of $[0] \in (\roman H_\phi^1)_{\roman{red}}$
onto an open neighborhood 
$U_\phi$
of
$[\phi] \in \roman{Rep}(\pi,G)$,
where
$W_\phi$ and
$U_\phi$
are endowed with the induced
stratified symplectic structures
$(C^{\infty}(W_\phi),\{\cdot,\cdot\})$
and
$(C^{\infty}(U_\phi),\{\cdot,\cdot\})$,
respectively (the notation
$\{\cdot,\cdot\}$ being abused somewhat).
\endproclaim

Thus locally, that is, 
near the point $[\phi]$ of $\roman{Rep}(\pi,G)$,
it suffices  to study
the reduced space
$(\roman H_\phi^1)_{\roman{red}}$
(written $\roman H_\phi$ in \cite\smooth).
The latter may
be understood by means of geometric invariant theory:
Write $W=\roman H_\phi^1$ and $K= Z_\phi$.
By a theorem of Hilbert and Hurwitz
\cite\weylbook\  (VIII \S 14),
the (graded) algebra $\Bobb R[W]^K$ of
$K$-invariant polynomials on $W$ 
has a finite set $p_1,\dots, p_d$ of generators.
Let 
$$
p=(p_1,\dots, p_d) \colon W @>>> \Bobb R^d
$$
be the Hilbert map, let $I$ be the ideal of relations among the $p_j$
in
$\Bobb R[y_1,\dots, y_d]$,
and let $Z$ be the corresponding algebraic subset of $\Bobb R^d$.
By the {\it Tarski-Seidenberg\/} Theorem,
$X:= \roman{Im}(p)$ is a semi-algebraic subset of $Z$,
cf. e.~g. \cite\bocosroy,
and the induced map
$\overline p \colon W\big/ K \to X$
is a homeomorphism.
With the appropriate structures, it is even a diffeomorphism
of smooth spaces \cite\gwschwar. 
If $S$ is a real affine $K$-subvariety of $W$,
it may be shown that 
there  is an algebraic subset $Z'$ of $Z$
such that the inequalities determining $X\subseteq Z$
determine $p(S)=S\big/ K:= X' =X \cap Z' \subseteq Z'$.
See e.~g. 
\cite\gwschwat\ for details.
Since the zero locus 
of $\Theta_{\phi}$ is 
a real affine $Z_\phi$-subvariety of $\roman H_\phi^1$,
this shows that
$\roman{Rep}(\pi,G)$ is locally semi-algebraic.

Near any of its points,
by means of the local model,  
we can now elucidate the stratified symplectic structure 
of $\roman{Rep}(\pi,G)$ (locally) as follows:
After identification of
$\roman H_\phi^*$
with the cohomology of $\Sigma$ with the appropriate coefficients,
and
after a choice of complex structure
on $\Sigma$ has been made,
the star operator endows $\roman H^1_\phi$ with a complex structure
compatible with the symplectic structure
$\sigma_\phi$, and in this way
$\roman H^1_\phi$
becomes a unitary $Z_\phi$-representation
such that
$\Theta_\phi$ is its
unique momentum mapping having the value zero at the origin.
In (2.4) on p.~52
of {\smc Arms-Gotay-Jennings}~\cite\armgojen,
this is called the \lq\lq the standard example\rq\rq.  
The upshot is that, locally,
the space
$\roman{Rep}(\pi,G)$
may be studied by techniques coming from
constrained systems in mechanics.

To simplify the exposition,
consider a general finite dimensional unitary representation
$W$ of a compact Lie group $K$.
Associated with it is the unique momentum mapping
$\mu$ from $W$ to $\fra k^*$
having the value zero at the origin \cite\armgojen.
In terms of complex coordinates $\bold z=(z_1,\dots,z_m)$
on $W$, 
when $K$ is a subgroup of $\roman U(m)$,
this momentum mapping is given by the formula
$$
(\xi \circ \mu) (\bold z) = \frac i2\sum \xi_{j,k} \overline z_j z_k
$$
where $\xi = \left(\xi_{j,k}\right) = \left(-\overline \xi_{k,j}\right)$
is a complex $(m \times m)$-matrix in $\fra k$;
here $\xi \in \fra k$ is viewed as a coordinate function on $\fra k^*$.

The $K$-action extends to an action 
of the complexification $K^{\Bobb C}$ of $K$ on $W$.
The idea is now to relate
the reduced space
$W_{\roman{red}} $
with the affine {\it categorical\/} quotient
$W \big / \big / K^{\Bobb C}$,
cf. e.~g. \cite\gwschwat.
Recall that, in general,
given a 
reductive group $\Gamma$ and an affine $\Gamma$-variety $V$,
the affine categorical quotient,
also referred to as the algebraic quotient
in the literature,
is a morphism $\tau\colon V \to V\big / \big /\Gamma$
of affine varieties which is constant on $\Gamma$-orbits and has the
following universal property:
Given a morphism $\psi \colon V \to V'$
of affine varieties which is constant on $\Gamma$-orbits,
there is a unique morphism
$\Psi \colon V\big / \big /\Gamma \to V'$
of affine varieties
such that $\Psi \circ \tau = \psi$.
With a slight abuse of language,
the variety $V\big / \big /\Gamma$
is then referred to as the categorical quotient as well.
Under our circumstances,
the 
categorical quotient
is 
the complex affine variety corresponding to
the algebra
$\Bobb C[W]^{K^{\Bobb C}}$
of 
$K^{\Bobb C}$-invariant complex polynomials;
actually 
the algebra of 
$K^{\Bobb C}$-invariant complex polynomials
coincides with that of $K$-invariant
complex polynomials.
By the  theorem of Hilbert and Hurwitz,
this algebra has a finite set $f_1,\dots,f_t$
of generators.
These yield 
a $K^{\Bobb C}$-invariant algebraic
map $f$ from $W$ to $\Bobb C^t$
which, by construction, factors through the space
of $K^{\Bobb C}$-orbits in $W$,
and the image $Y$
in $\Bobb C^t$ is the variety defined by the relations
among the
$f_1,\dots,f_t$.
This variety is a model for $W \big / \big / K^{\Bobb C}$.
See e.~g. \cite\gwschwat\ (\S 3) for details.
An observation of {\smc Kempf-Ness} \cite\kemneone,
cf. \S 4 of \cite\gwschwat, where 
the zero locus $\mu^{-1}(0)$ is referred to as a {\it Kempf-Ness\/}
set,
implies that
the canonical map
$
W_{\roman{red}} @>>>
W\big /\big/K^{\Bobb C}
$
from
the reduced space
$W_{\roman{red}}$
to the (affine) categorical quotient
$W\big /\big/K^{\Bobb C}$
induced by the inclusion of
$\mu^{-1}(0)$ into $W$
is a homeomorphism. 
As a space, in fact, as a complex affine variety,
the reduced space
$W_{\roman{red}}$
thus looks like the affine categorical quotient
$W\big /\big/K^{\Bobb C}$. As a stratified symplectic space
it looks somewhat different, though.
We mention in passing that the stratified symplectic and the K\"ahler 
structures combine to what we called a stratified K\"ahler structure in
\cite\kaehler.

We now illustrate
this local model analysis
for the special case
where $G=\roman{SU}(2)$.
We remind the reader that the decomposition
$$
N = N_G \cup N_{(T)}\cup N_Z
$$
has been described earlier.

\proclaim{Theorem 4}
Near a point of $N_{(K)}$,
$N$ and
$(C^{\infty}(N),\{\cdot,\cdot\})$
may be described in the following way:
\smallskip
\noindent
$K=Z$: the space $\Bobb C^{3(\ell -1)}$ with its standard symplectic
Poisson structure;
\newline\noindent
$K=T$: a product of 
$\Bobb C^{\ell}$ with its standard symplectic
Poisson structure
and of the reduced space and reduced Poisson algebra
of a system of $\ell -1$ particles in the plane
with total angular momentum zero;
\newline\noindent
$K=G$: the reduced space and reduced Poisson algebra
of a system of $\ell$ particles in 3-space
with total angular momentum zero.
\endproclaim

The proof of this theorem 
relies on the local models:
For example,
when $\ell = 2$, 
near a point of the \lq\lq middle\rq\rq\  stratum $N_{(T)}$,
the space $N$ looks like the product of a copy of
$\Bobb C^2$ with the reduced system of a single particle
in the plane $\Bobb R^2$.
The latter is $\Bobb R^2$,
endowed with the smooth structure
and Poisson algebra (2.2).

In some more detail:
For a point $[\phi]$ in the top stratum,
$\roman H^0_\phi$ and 
$\roman H^2_\phi$
are zero,
and hence near $[\phi]$, the moduli space looks like
a neighborhood of zero in $\roman H^1_\phi$,
with the symplectic structure $\sigma_\phi$.
The latter, in turn,
amounts
to a copy of
$\Bobb C^{3(\ell-1)}$, with the standard symplectic structure.

For a point $[\phi]$ in the middle  stratum
$N_{(T)}$,
$\fra g_{\phi}$ 
decomposes into a direct sum
of $\fra t$ and $\fra t^{\bot}$ where
$\fra t$ is the Lie algebra
of $T$
which is a copy 
of the reals 
(with trivial $T$-action)
and
$\fra t^{\bot}$
amounts to
$\Bobb R^2$,
with circle action through the 2-fold covering
map onto $\roman{SO}(2,\Bobb R)$.
Moreover, $\roman H_\phi^1$
decomposes into a direct sum of
$(\fra t \otimes \Bobb C)^\ell$
and
$(\fra t^{\bot} \otimes \Bobb C)^{\ell-1}$.
The action on 
$(\fra t \otimes \Bobb C)^\ell$
is trivial -- in fact, this summand corresponds
to the points {\it in\/} the stratum $N_{(T)}$
(locally),
while the
$\roman{SO}(2,\Bobb R)$-representation
on
$(\fra t^{\bot} \otimes \Bobb C)^{\ell-1}$
is hamiltonian, with momentum mapping
$\Theta_\phi$, restricted to
$(\fra t^{\bot} \otimes \Bobb C)^{\ell-1}$.
The latter boils down
to the classical constrained system
of $\ell -1$
particles moving in the plane with constant total
angular momentum.
In particular,
reduction at total angular momentum zero
yields the reduced space we are looking for
(locally).

When $\ell = 2$,
$W = \fra t^{\bot} \otimes \Bobb C \cong \Bobb R^2 \times \Bobb R^2$.
The 
$\roman{SO}(2,\Bobb R)$-representation
is the obvious one,
that is,
$\roman{SO}(2,\Bobb R)$
acts as rotation group on each copy of
$\Bobb R^2$.
With the usual coordinates
$(q,p) \in \Bobb R^2 \times \Bobb R^2$,
the momentum mapping 
$\mu$ from $W$ to $\Bobb R$
is given by the assignment
to $(q,p)$ of
the determinant $|q p |$.
The algebra of (real) invariants 
in $\Bobb R[W]$
is generated 
by the three scalar products
$qq$, $qp$, $pp$, and the determinant
$|qp|$.
However, on the zero locus
$\mu^{-1}(0)$,
the determinant vanishes 
whence 
the algebra 
$C^{\infty}(W_{\roman{red}})$
is generated by the 
three scalar products.

To illustrate how we can understand
$W_{\roman{red}}$ as a {\it space\/}
by means of the corresponding  categorical
affine quotient
we observe that the extension of
the 
$\roman{SO}(2,\Bobb R)$-representation
to its complexification
amounts to the standard
$\roman{SO}(2,\Bobb C)$-representation
on $\Bobb C^2$.
This representation has a single
invariant,
the complex scalar product 
which, with the notation $w = q + ip$, we write
$ww$,
and the algebra of complex invariants is free.
Consequently the affine categorical quotient
$\Bobb C^2 \big / \big /\roman{SO}(2,\Bobb C)$
is a copy of the complex line $\Bobb C$. 
By virtue of the observation of
Kempf-Ness,
the
canonical map
from
$W_{\roman{red}}$
to $\Bobb C$
is a homeomorphism.
Under this homeomorphism,
with the notation $ww = x_1 + i x_2$,
 so that $x_1$ and $x_2$ are the coordinate
functions on $\Bobb C$,
viewed as the real plane,
we have
$$
x_1 = qq - pp,
\ 
x_2 = 2 qp,
\ 
\rho = qq + pp.
$$
It is obvious that the algebra 
$C^{\infty}(W_{\roman{red}})$
is as well generated by the coordinate
functions $x_1 ,x_2$ and the radius function $\rho$.
Moreover, the complex picture tells us that
the single obvious relation (2.1)
between the coordinate functions
and the radius function 
is defining,
that is, the relation
$x_1^2 + x_2^2 = \rho^2$ suffices;
hence the latter is a defining relation for
$C^{\infty}(W_{\roman{red}})$.
Finally, a straightforward calculation
of the Poisson brackets between
$x_1, x_2, \rho$,
viewed as functions
on the {\it original\/} space $W$,
yields the formulas
already spelled out as (2.2).
See \cite\locpois\ for details.
Since the Zariski tangent space
for the non-standard structure 
$C^{\infty}(W_{\roman{red}})$
on $\Bobb R^2= W_{\roman{red}}$
is a copy of $\Bobb R^3$,
the Zariski tangent space
at a point 
$[\phi]$ of the middle stratum
$N_{(T)}$
is a copy of $\Bobb R^7$. 

A similar reasoning yields
a model for a neighborhood
of a point in the \lq\lq bottom\rq\rq\ 
stratum $N_G$.
This stratum consists of
$2^{2\ell}$
isolated points
and, locally, the Poisson algebra
is that of the reduced classical constrained
system of $\ell$ particles
in $\Bobb R^3$ with total angular momentum zero.
Even for $\ell =2$,
the reduced
Poisson algebra is already rather
complicated:
It has ten generators;
in fact, 
the elements of a basis of
$\fra{sp}(2,\Bobb R)$
may be taken as coordinate functions
on $\fra{sp}(2,\Bobb R)^*$,
and the reduced space may be described
as the closure of a certain nilpotent orbit
in
$\fra{sp}(2,\Bobb R)^*$.
This relies on the theory
of dual pairs \cite\howeone\  and is worked out in
\cite\lermonsj;
see our paper \cite\locpois\  for details.
A more general theory explaining how 
and why nilpotent orbits come into play here
has been developed in \cite\kaehler.
Suffice it to mention at this stage that,
once
$\fra{sp}(2,\Bobb R)$ has been identified with its dual
$\fra{sp}(2,\Bobb R)^*$ by means of the Killing form,
the
$\roman{Sp}(2,\Bobb R)$-momentum mapping
$$
(\roman T^* \Bobb R^3)^{\times 2}
\cong
\Bobb R^3 \otimes \Bobb R^4
@>>>
\fra{sp}(2,\Bobb R)^*
$$
for the induced hamiltonian
$\roman{Sp}(2,\Bobb R)$-action 
on $(\roman T^* \Bobb R^3)^{\times 2}$
coming from the standard representation
of $\roman{Sp}(2,\Bobb R)$ 
on $\Bobb R^4$ 
identifies the
$\roman{O}(3,\Bobb R)$-reduced space,
that is, the reduced space of two particles in
$\Bobb R^3$
with total angular momentum zero,
with the closure of the nilpotent orbit generated by
$$
\left [
\matrix
0&0&1&0
\\
0&0&0&1
\\
0&0&0&0
\\
0&0&0&0
\endmatrix
\right]
\in 
\fra{sp}(2,\Bobb R).
$$
The Zariski tangent space
at any of the 16 double points of the Kummer surface,
for the {\it non-standard\/} smooth structure on
$\roman{Rep}(\pi, \roman{SU}(2))$
we are presently discussing,
has dimension 10;
see what is said in the next section.
Topologically, these 16 points are not singularities, though:
as a space,
$\roman{Rep}(\pi, \roman{SU}(2))$ is just
complex projective 3-space
which, with its {\it ordinary\/}
smooth structure,
at any of these points,
is plainly non-singular.

\medskip\noindent
{\bf 6. Cleaning up the Zariski tangent space}
\smallskip\noindent
The notions of singularity and Zariski tangent space
refer to a particular structure.
The representation spaces
described above carry real (locally semi-algebraic) structures
as well as complex analytic ones;
at a point of such a space, the 
Zariski tangent space
for 
the complex structure will
in general {\it not\/} be 
the 
complexification of
the Zariski tangent space
for the real structure,
and 
the lack of distinction between
various structures
seems to have created some confusion in the literature.
We take the opportunity to try to \lq\lq clean up\rq\rq\ 
the situation.

As before, let $G$ be a 
compact Lie group and $\pi$ 
the fundamental group 
of a closed surface $\Sigma$. 
In 
\cite\goldmone\ (p.~205),
\cite\kapmione\ (p.~1091),
and elsewhere in the literature,
for a homomorphism $\phi \colon \pi \to G$
(where $G$ could be any real algebraic group and
$\pi$ any finitely presented discrete group),
the vector space $\roman H^1(\pi,\fra g_{\phi})$
is referred to as 
the {\it Zariski tangent space\/} of
$\roman{Rep}(\pi,G)$ at $[\phi]$;
in 
\cite\goldmone,
\lq Zariski tangent space\rq\ 
is actually put in quotation marks.
In \cite\brymcone\ (proof of Lemma 11 on p.~49),
it is asserted that, given 
any homomorphism $\phi$
from  
$\pi$ 
to $G$,
the vector space $\roman H^1(\pi,\fra g_{\phi})$
is the \lq tangent space\rq\  of
$\roman{Rep}(\pi,G)$ at $[\phi]$
(whether or not $\roman{Rep}(\pi,G)$ is smooth near $[\phi]$).
None of the cited references makes precise, with reference
to which structure, i.~e. algebra of functions, 
the Zariski tangent space is to be taken.
What seems intended 
is this:
Write  $\ell$
for the genus of $\Sigma$, let 
$\phi \in \roman {Hom}(\pi, G) \subseteq G^{2 \ell}$,
and consider the derivative $dr_{\phi}$
of the word map $r\colon G^{2 \ell} @>>> G$
at $\phi$; right translation identifies
the kernel $\roman{ker} (dr_{\phi})$
of $dr_{\phi}$
with the space
of 1-cocycles for $\pi$ in $\fra g_{\phi}$
\cite\goldmone\ and, by Theorem 7.14 of \cite\smooth,
$\roman{ker} (dr_{\phi})$ amounts to the Zariski tangent space
of 
$\roman {Hom}(\pi, G)$
with reference to the smooth structure
induced from the embedding into $G^{2 \ell}$.
Moreover,
right translation identifies
the tangent space 
$\roman T_{\phi}G\phi$
at $\phi$
to the $G$-orbit $G\phi$ in
$\roman {Hom}(\pi, G) \subseteq G^{2 \ell}$
with the space of 1-coboundaries
for $\pi$ in $\fra g_{\phi}$; in the  
cited references,
$\roman H^1(\pi,\fra g_{\phi})$
is apparently viewed as the normal space
to $\roman T_{\phi}G\phi$
inside 
$\roman{ker} (dr_{\phi})$,
and the Zariski tangent space
of
$\roman{Rep}(\pi,G)$ at $[\phi]$
seems to have been confused with this normal space.
At a point $[\phi]$ of the top stratum, the notion of normal space
indeed makes sense and 
$\roman H^1(\pi,\fra g_{\phi})$
comes down to the ordinary tangent space,
in particular, may be identified with the
normal space
to $\roman T_{\phi}G\phi$.

However,
the first objection to taking $\roman H^1(\pi,\fra g_{\phi})$
as the Zariski tangent space
in general
is the purely formal observation that
$\roman H^1(\pi,\fra g_{\phi})$
is not even invariantly defined in terms of
the point $[\phi]$ of
$\roman{Rep}(\pi,G)$
unless 
$[\phi]$ belongs to the top stratum of
$\roman{Rep}(\pi,G)$.
More precisely,
given the two representatives
$\phi$ and $\phi'=\roman{Ad}(x)\phi$
of $[\phi]$, where $x \in G$,
the association $w \mapsto \roman{Ad}(x)w$ ($w \in \fra g$)
yields an isomorphism
$x^{\sharp}\colon\fra g_{\phi} \to \fra g_{\phi'}$
of $\pi$-modules but this isomorphism depends on the choice
of $x$; in particular, when
$\phi'=\phi$, that is,
when $x$ lies in the stabilizer $Z_{\phi}\subseteq G$ of $\phi$,
there is no need for the
automorphism
$x^{\sharp}\colon\fra g_{\phi} \to \fra g_{\phi}$
to be the identity;
indeed, the resulting $Z_{\phi}$-representation
on $\fra g$ ($=\fra g_{\phi}$ but 
the fact that
$\fra g$
carries the $\pi$-representation $\phi$ is not relevant at this point)
will in general be non-trivial.
For example, under the circumstances of Theorem 4,
the \lq\lq bottom\rq\rq\ stratum consists of
$2^{2\ell}$
isolated points, and the stabilizer
$Z_{\phi}$ for any representative $\phi$
of any of these points coincides with the whole group $G=\roman{SU}(2)$,
and the adjoint representation of $G$ on $\fra g$
is plainly non-trivial.
Likewise, under these circumstances,
the stabilizer
$Z_{\phi}$ for any representative $\phi$
of a point $[\phi]$ of the middle  stratum
is a maximal torus $T$ in 
$\roman{SU}(2)$,
and the resulting $T$-representation
is plainly non-trivial.

We now explain briefly the relationship between
$\roman H^1(\pi,\fra g_{\phi})$
and the Zariski tangent space
$\roman T_{[\phi]}\roman{Rep}(\pi,G)$.
More details and proofs may be found in Section 7 of \cite\smooth.
Theorem 3
reduces the
relationship under discussion
to that between
$\roman H^1(\pi,\fra g_{\phi})$ and
the Zariski tangent space
$\roman T_{[0]}(\roman H_\phi^1)_{\roman{red}}$.
To elucidate the latter relationship, let $V_{\phi}\subseteq \roman H_\phi^1$ 
be the zero locus of the momentum mapping
$\Theta_{\phi}$.
The projection from
$V_{\phi}$
to
$(\roman H_\phi^1)_{\roman{red}}$
induces a linear map
$\roman T_0V_{\phi} \to \roman T_{[0]}(\roman H_\phi^1)_{\roman{red}}$
between the Zariski tangent spaces.
Furthermore, by Lemma 7.6 in \cite\smooth, 
$V_{\phi}$
spans
$\roman H_\phi^1$ whence
$\roman T_0V_{\phi} =
\roman H_\phi^1$; thus
the projection from
$V_{\phi}$
to
$(\roman H_\phi^1)_{\roman{red}}$
induces a linear map
$$
\lambda
\colon
\roman H_\phi^1 @>>> \roman T_{[0]}(\roman H_\phi^1)_{\roman{red}}.
$$
By Theorem 3,
the choice
of $\phi$
(in its class $[\phi]$)
induces a diffeomorphism
of an open neighborhood
$W_\phi$ of $[0] \in (\roman H_\phi^1)_{\roman{red}}$
onto an open neighborhood 
$U_\phi$
of
$[\phi] \in \roman{Rep}(\pi,G)$
and hence
an isomorphism from
$\roman T_{[0]}(\roman H_\phi^1)_{\roman{red}}$
onto
$\roman T_{[\phi]}\roman{Rep}(\pi,G)$;
the latter combines with $\lambda$ to 
a linear map
$$
\lambda_{\phi}
\colon
\roman H_\phi^1 @>>> 
\roman T_{[\phi]}\roman{Rep}(\pi,G).
$$
The map $\lambda_\phi$ has the following properties
(see p. 214 of \cite\smooth):
\newline\noindent {\rm (1)}
{\sl It is independent of the choice of
$\phi$ in the sense that,
for every $x\in G$, the 
composite
$$
\roman H^1(\pi,g_{\phi})
@>{\roman{Ad}_\flat(x)}>>
\roman H^1(\pi,g_{x\phi})
@>{\lambda_{x\phi}}>>
\roman T_{[\phi]} \roman {Rep}(\pi,G)
$$
of the induced linear isomorphism $\roman{Ad}_\flat(x)$
with 
$\lambda_{x\phi}$ 
coincides with
$\lambda_{\phi}$.
\newline\noindent {\rm (2)}
Its kernel equals the subspace
of $\roman H^1(\pi,g_{\phi})$
which is the linear span of the elements $xw-w$
where $w \in \roman H^1(\pi,g_{\phi})$ and
$x \in Z_{\phi}$.
\newline\noindent {\rm (3)}
Its image equals the 
(smooth)
tangent space
$\roman T_{[{\phi}]}(\roman {Rep}(\pi,G)_{(K)})$, viewed as a subspace
of
$\roman T_{[{\phi}]}\roman {Rep}(\pi,G)$
in an obvious sense,
where
$\roman {Rep}(\pi,G)_{(K)}$ denotes the stratum in which $[\phi]$ lies,
so that $K= Z_{\phi}$.
\newline\noindent {\rm (4)}
It is an isomorphism
if and only if $[{\phi}]$ is a non-singular point of
$\roman {Rep}(\pi,G)$, i.~e. belongs to the top stratum.}

Under the circumstances of Theorem 4,
for the special case where the genus
$\ell$ of the surface $\Sigma$
equals  2,
these observations entail the following insight
into the
Zariski tangent spaces of the
corresponding space $N=\roman {Rep}(\pi,G)$:
For
a point $[\phi]$ of the middle stratum
$N_{(T)}$,
the Zariski tangent space
$\roman T_{[\phi]} N$ has dimension 4 + 3 = 7.
On the other hand,
the dimension of $\roman H_\phi^1$
equals 8, and
the linear map
$\lambda_\phi$ from
$\roman H_\phi^1$
to
$\roman T_{[\phi]} N$
has rank four.
Thus 
the Zariski tangent space
$\roman T_{[\phi]} N$
can
in no way be identified with
the cohomology group 
$\roman H_\phi^1$.
Likewise,
let
$[\phi]$  be a point 
in 
$N_{G}$.
Then the Zariski tangent space
$\roman T_{[\phi]} N$ has dimension 10 and hence
the minimal number of generators
of
$C^{\infty}(N)$
near $[\phi]$ or rather that of its germ at $[\phi]$ is 10;
see p.~217 of \cite\smooth\ for details.
Moreover, a closer look reveals that
the Zariski tangent space
$\roman T_{[\phi]} N$
equals that
of
$\roman T_{[\phi]} N_{(T)}$, with reference to
the induced smooth structure
$C^{\infty}(N_{(T)})$.
In fact,
in the language of constrained systems,
in the local model, $N_{(T)}$
corresponds to reduced states
where each of the two particles individually
has angular momentum zero,
cf. what is said in our paper
\cite\locpois.
In particular,
the minimal number of generators
of the induced smooth structure
$C^{\infty}(N_{(T)})$
near $[\phi]$ or rather that of its germ at $[\phi]$ is still 10.
Finally,
the linear map
$\lambda_\phi$ from
$\roman H_\phi^1$
to
$\roman T_{[\phi]} N$
is zero
since the derivative of $\lambda$
at the origin
is zero.
Thus 
the Zariski tangent space
$\roman T_{[\phi]} N$
can
in no way be identified with
the cohomology group 
$\roman H_\phi^1$, which has dimension 12.
As for the complex analytic structure,
we recall that, cf. \cite\locpois,
as a complex variety,
near a point
$[\phi]$ 
in 
$N_{G}$,
the stratum
$N_{(T)}$
looks like the quadric
$Y^2 = XZ$
in complex 3-space.
Hence,
at a point
$[\phi]$ 
in 
$N_{G}$,
the complex Zariski
tangent space
of
$N_{(T)}$
has dimension 3,
and this 
Zariski
tangent space
coincides with the smooth complex tangent space
of $N$ at $[\phi]$.
But,
as noted above, the real
Zariski tangent space $\roman T_{[\phi]} N$
at the point $[\phi]$
has dimension 7,
and $\roman T_{[\phi]} N$
actually coincides with
$\roman T_{[\phi]} N_{(T)}$.
Thus, as a smooth space,
with smooth structure
$C^{\infty}(\roman{Rep}(\pi,G))$
given in Theorem 1,
the space $N=\roman {Rep}(\pi,G)$ 
looks rather different
from complex projective 3-space
with its standard smooth structure.

The observation that the tangent cones
for varieties of spaces
of homomorphisms $\roman{Hom}(\pi,G)$
for suitable discrete groups $\pi$
(e.~g. fundamental groups of compact K\"ahler manifolds)
and appropriate Lie groups $G$
are quadratic as well
as the rigidity results for such spaces,
due to {\smc Goldman and Millson\/} and others,
see \cite\golmilon\ 
and the references there,
were influential in the development of
the subject,
and these results are unaffected by our remarks.
In particular, the vector space
$\roman H_\phi^1$
or what corresponds to it 
under certain more general circumstances
is
a constituent of a differential graded
Lie algebra controlling the corresponding deformation problem
under consideration; see e.~g. \cite{\goldmone,\,\golmiltw,\,\kapmione}.

We conclude with a comment on another abuse 
of language:
In the literature,
it is common to refer to a space
of the kind
$\roman{Rep}(\pi, G)$,
$G$ being a real algebraic Lie group,
as a \lq\lq representation variety\rq\rq. 
As explained above,
$\roman{Rep}(\pi, G)$ 
is {\it not\/} a real variety in an obvious way;
in particular,
$\roman{Rep}(\pi, G)$ 
does {\it not\/} consist of the real points of the
corresponding complex representation variety
$\roman{Hom}(\pi, G_{\Bobb C})\big / \big / G_{\Bobb C}$ 
(where
the notation $\big / \big /$ refers to the corresponding
categorical quotient).
The \lq\lq naive\rq\rq\ representation space
$\roman{Rep}(\pi, G) = \roman{Hom}(\pi, G)\big / G$
lies inside the
real categorical quotient
$\roman{Hom}(\pi, G)\big / \big / G$,
and the latter is indeed a real variety.
Locally, the difference between
$\roman{Rep}(\pi, G)$
and
$\roman{Hom}(\pi, G)\big / \big / G$
is of the kind as that between a half cone and a
(double) cone explained in Section 2 above.

\medskip
\centerline{References}

\medskip
\widestnumber\key{999}

\ref \no  \abramars
\by R. Abraham and J. E. Marsden
\book Foundations of Mechanics
\publ Benjamin-\linebreak
Cummings Publishing Company
\yr 1978
\endref

\ref \no  \armcusgo
\by J. M. Arms,  R. Cushman, and M. J. Gotay
\paper  A universal reduction procedure for Hamiltonian group actions
\paperinfo in: The geometry of Hamiltonian systems, T. Ratiu, ed.
\jour MSRI Publ. 
\vol 20
\pages 33--51
\yr 1991
\publ Springer
\publaddr Berlin $\cdot$ Heidelberg $\cdot$ New York $\cdot$ Tokyo 
\endref

\ref \no  \armgojen
\by J. M. Arms, M. J. Gotay, and G. Jennings
\paper  Geometric and algebraic reduction for singular
momentum mappings
\jour Advances in Mathematics
\vol 79
\yr 1990
\pages  43--103
\endref

\ref \no  \arnobook
\by V. I. Arnold
\book Mathematical Methods of classical mechanics
\bookinfo Graduate Texts in Mathematics, No. 60
\publ Springer
\publaddr Berlin $\cdot$ Heidelberg $\cdot$ New York $\cdot$ Tokyo 
\yr 1978, 1989 (2nd edition)
\endref

\ref \no  \atibottw
\by M. Atiyah and R. Bott
\paper The Yang-Mills equations over Riemann surfaces
\jour Phil. Trans. R. Soc. London  A
\vol 308
\yr 1982
\pages  523--615
\endref

\ref \no \berghueb
\by C. Berger and J. Huebschmann
\paper Comparison of the geometric bar and $W$-\linebreak
constructions
\jour J. of Pure and Applied Algebra
\vol 131
\yr 1998
\pages 109--123
\endref

\ref \no \bocosroy
\by J. Bochnak, M. Coste, and M.-F. Roy
\book G\'eom\'etrie
alg\'ebrique r\'eelle
\bookinfo Ergebnisse der Mathematik und ihrer Grenzgebiete, 3. Folge, Band 12.
A series of modern surveys in Mathematics
\vol 12
\yr 1987
\publ Springer
\publaddr Berlin $\cdot$ Heidelberg $\cdot$ New York  $\cdot$ Tokyo
\endref

\ref \no \brymcone
\by J. L. Brylinski and D. A. McLaughlin
\paper  Holomorphic quantization and unitary representations
\paperinfo in: Lie Theory and Geometry, 
In honor of B. Kostant,
J. L. Brylinski, R. Brylinski, V. Guillemin, V. Kac, eds. 
\jour Progress in Mathematics 
\vol 1994
\publ Birkh\"auser
\publaddr Boston $\cdot$ Basel $\cdot$ Berlin
\yr 1994
\pages 21--64
\endref

\ref \no \goldmone
\by W. M. Goldman
\paper The symplectic nature of the fundamental group of surfaces
\jour Advances
\vol 54
\yr 1984
\pages 200--225
\endref

\ref \no  \goldmtwo
\by W. M. Goldman
\paper Invariant functions on Lie groups and Hamiltonian flows of
surface group representations
\jour Inventiones
\vol 85
\yr 1986
\pages 263--302
\endref

\ref \no \golmilon
\by W. M. Goldman and J. Millson
\paper The deformation theory of representations of
fundamental groups of compact Kaehler manifolds
\jour Publ. Math. I. H. E. S.
\vol 67
\yr 1988
\pages 43--96
\endref

\ref \no \golmiltw
\by W. M. Goldman and J. Millson
\paper Differential graded Lie algebras and singularities
of level set momentum mappings
\jour Commun. Math. Phys. 
\vol 131
\yr 1990
\pages 495--515
\endref

\ref \no \gormacon
\by M. Goresky and R. MacPherson
\paper Intersection homology theory
\jour Topology
\vol 19
\yr 1980
\pages 135--162
\endref

\ref \no \guhujewe
\by K. Guruprasad, J. Huebschmann, L. Jeffrey, and A. Weinstein
\paper Group systems, groupoids, and moduli spaces
of parabolic bundles
\jour Duke Math. J.
\vol 89
\yr 1997
\pages 377--412
\endref

\ref \no \hartsboo
\by  R. Hartshorne
\book Algebraic Geometry
\bookinfo Graduate texts in Mathematics
 No. 52
\publ Springer
\publaddr Berlin-G\"ottingen-Heidelberg
\yr 1977
\endref

\ref \no  \howeone
\by R. Howe
\paper Remarks on classical invariant theory
\jour  Trans. Amer. Math. Soc.
\vol 313
\yr 1989
\pages  539--570
\endref

\ref \no \poiscoho
\by J. Huebschmann
\paper Poisson cohomology and quantization
\jour 
J. f\"ur die Reine und Angewandte Mathematik
\vol  408 
\yr 1990
\pages 57--113
\endref

\ref \no  \souriau
\by J. Huebschmann
\paper On the quantization of Poisson algebras
\paperinfo in:
Symplectic Geometry and Mathematical Physics,
Actes du colloque en l'honneur de Jean-Marie Souriau,
P. Donato, C. Duval, J. Elhadad, G.M. Tuynman, eds. 
\jour Progress in Mathematics 
\vol 99
\publ Birkh\"auser
\publaddr Boston $\cdot$ Basel $\cdot$ Berlin
\yr 1991
\pages 204--233
\endref

\ref \no \singula
\by J. Huebschmann
\paper The singularities of Yang-Mills connections
for bundles on a surface. I. The local model
\jour Math. Z. 
\vol 220
\yr 1995
\pages 595--605
\endref

\ref \no \singulat
\by J. Huebschmann
\paper The singularities of Yang-Mills connections
for bundles on a surface. II. The stratification
\jour Math. Z. 
\vol 221
\yr 1996
\pages 83--92
\endref

\ref \no \smooth
\by J. Huebschmann
\paper 
Smooth structures on moduli spaces of central Yang-Mills connections 
for bundles on a surface
\jour J. of Pure and Applied Algebra
\vol 126
\yr 1998
\pages 183--221
\endref

\ref \no \poisson
\by J. Huebschmann
\paper 
Poisson
structures on certain
moduli spaces 
for bundles on a surface
\jour Annales de l'Institut Fourier
\vol 45
\yr 1995
\pages 65--91
\endref

\ref \no \locpois
\by J. Huebschmann
\paper Poisson geometry of flat connections 
for {\rm SU(2)}-bundles on surfaces
\jour Math. Z.
\vol 221
\yr 1996
\pages 243--259
\endref

\ref \no \modus
\by J. Huebschmann
\paper Symplectic and Poisson structures of certain moduli spaces
\jour Duke Math. J.
\vol 80
\yr 1995
\pages 737--756
\endref

\ref \no \modustwo
\by J. Huebschmann
\paper Symplectic and Poisson structures of certain moduli spaces. II.
Projective
representations of cocompact planar discrete groups
\jour Duke Math. J.
\vol 80
\yr 1995
\pages 757--770
\endref

\ref \no \moscou
\by J. Huebschmann
\paper
Poisson geometry of certain
moduli spaces for bundles on a surface
\paperinfo
\lq\lq A translation of algebra\rq\rq,
International Geometrical Colloquium,
Moskau, 1993;
Vseross. Inst. Nauchn. i. Tekhn. Inform. (VINITI)
Moscou 1995,
edited by E. I. Kuznetsova
\jour J. Math. Sci. 
\vol 82 
\yr 1996
\pages 3780--3784
\endref

\ref \no \srni
\by J. Huebschmann
\paper
Poisson geometry of certain
moduli spaces
\paperinfo
Lectures delivered at the \lq\lq 14th Winter School\rq\rq, Srni,
Czeque Republic,
January 1994
\jour Rendiconti del Circolo Matematico di Palermo, Serie II
\vol 39
\yr 1996
\pages 15--35
\endref

\ref \no \claustha
\by J. Huebschmann
\paper On the Poisson geometry of certain moduli spaces
\paperinfo in: Proceedings of an international workshop on
\lq\lq Lie theory and its applications in physics\rq\rq,
Clausthal, 1995
H. D. Doebner, V. K. Dobrev, J. Hilgert, eds.
\publ World Scientific
\publaddr Singapore $\cdot$
New Jersey $\cdot$
London $\cdot$
Hong Kong 
\pages 89--101
\yr 1996
\endref

\ref \no \kan
\by J. Huebschmann
\paper Extended moduli spaces, the 
Kan construction, and lattice gauge theory
\jour Topology
\vol 38
\yr 1999
\pages 555--596
\endref

\ref \no \poismod
\by J. Huebschmann
\paper On the variation of the Poisson structures of certain moduli spaces
\paperinfo dg-ga/9710033 
\jour Math. Ann. (to appear)
\endref

\ref \no \kaehler
\by J. Huebschmann
\paper Kaehler spaces, nilpotent orbits, and singular reduction
\paperinfo in preparation
\jour 
\endref

\ref \no \huebjeff
\by J. Huebschmann and L. Jeffrey
\paper Group Cohomology Construction of Symplectic Forms
on Certain Moduli Spaces
\jour Int. Math. Research Notices
\vol 6
\yr 1994
\pages 245--249
\endref

\ref \no \kapmione
\by M. Kapovich and J. J. Millson
\paper On the deformation theory of representations of fundamental groups
of compact hyperbolic 3-manifolds
\jour Topology
\vol 35
\yr 1996
\pages 1085--1106
\endref

\ref \no \kemneone
\by G. Kempf and L. Ness
\paper The length of vectors in representation spaces
\jour Lecture Notes in Mathematics
\vol 732
\yr 1978
\pages 233--244
\paperinfo Algebraic geometry, Copenhagen, 1978
\publ Springer 
\publaddr Berlin $\cdot$ Heidelberg $\cdot$ New York
\endref

\ref \no \kobanomi
\by S. Kobayashi and K. Nomizu
\book Foundations of differential geometry, I (1963), II (1969)
\bookinfo Interscience Tracts in Pure and Applied Mathematics, No. 15
\publ Interscience Publ.
\publaddr New York-London-Sydney
\endref

\ref \no \kunzbotw
\by E. Kunz
\book Einf\"uhrung in die kommutative Algebra und algebraische Geometrie
\bookinfo Vieweg Studium Band 46 Aufbaukurs Mathematik
\publ Friedrich Vieweg \& Sohn
\publaddr Braunschweig/Wiesbaden
\yr 1980
\moreref
\book engl. translation: Introduction
to commutative algebra and algebraic geometry
\publ Birkh\"auser
\publaddr Boston
\yr 1985
\endref

\ref \no \lermonsj
\by E. Lerman, R. Montgomery and R. Sjamaar
\paper Examples of singular reduction
\book Symplectic Geometry 
\bookinfo Warwick, 1990, ed. D. A. Salamon, London Math. Soc. Lecture Note 
Series
\vol 192
\yr 1993
\pages  127--155
\publ Cambridge University Press
\publaddr Cambridge, UK
\endref

\ref \no \marswein
\by J. Marsden and A. Weinstein
\paper Reduction of symplectic manifolds with symmetries
\jour Rep. on Math. Phys.
\vol 5
\yr 1974
\pages 121--130
\endref

\ref \no \naramntw
\by M. S. Narasimhan and S. Ramanan
\paper Moduli of vector bundles on a compact Riemann surface
\jour Ann. of Math.
\vol 89
\yr 1969
\pages  19--51
\endref

\ref \no \naramnth
\by M. S. Narasimhan and S. Ramanan
\paper 2$\theta$-linear systems on abelian varieties
\jour Bombay colloquium
\yr 1985
\pages  415--427
\endref

\ref \no \narashed
\by M. S. Narasimhan and C. S. Seshadri
\paper Stable and unitary vector bundles on a compact Riemann surface
\jour Ann. of Math.
\vol 82
\yr 1965
\pages  540--567
\endref

\ref \no \newstboo
\by P. E. Newstead
\book Introduction to Moduli Problems and Orbit Spaces
\bookinfo Tata Institute Lecture Notes
\publ Springer
\publaddr Berlin--G\"ottingen--Heidelberg
\yr 1978
\endref

\ref \no \gwschwar
\by G.W. Schwarz
\paper Smooth functions invariant under the action of
a compact Lie group
\jour Topology 
\vol 14
\yr 1975
\pages 63--68
\endref

\ref \no \gwschwat
\by G. W. Schwarz
\paper The topology of algebraic quotients
\paperinfo In: Topological methods in algebraic
transformation groups
\jour Progress in Math.
\vol 80
\yr 1989
\pages 135--152
\publ Birkh\"auser
\publaddr Boston $\cdot$ Basel $\cdot$ Berlin
\endref

\ref \no \seshaone
\by C. S. Seshadri
\paper Spaces of unitary vector bundles on a compact Riemann surface
\jour Ann. of Math.
\vol 85
\yr 1967
\pages 303--336
\endref

\ref \no \seshaboo
\by C. S. Seshadri
\book Fibr\'es vectoriels sur les courbes alg\'ebriques
\bookinfo Ast\'erisque Vol. 96
\publ Soc. Math. de France
\yr 1982
\endref

\ref \no \sjamlerm
\by R. Sjamaar and E. Lerman
\paper Stratified symplectic spaces and reduction
\jour Ann. of Math.
\vol 134
\yr 1991
\pages 375--422
\endref

\ref \no  \weiltwo
\by A. Weil
\paper Remarks on the cohomology of groups
\jour                                      
Ann. of Math.
\vol 80                                    
\yr 1964
\pages  149--157
\endref

\ref \no  \weinstwo
\by A. Weinstein
\paper The local structure of Poisson manifolds
\jour J. of Diff. Geom.
\vol 18
\yr 1983
\pages 523--557
\endref

\ref \no  \weinsone
\by A. Weinstein
\paper Poisson structures
\jour Ast\'erisque,
\vol hors-serie, 
\yr 1985
\pages 421--434
\paperinfo in E. Cartan et les Math\'ematiciens d'aujourd hui, 
Lyon, 25--29 Juin, 1984
\endref

\ref \no  \weinsthi
\by A. Weinstein
\paper The symplectic structure on moduli space
\paperinfo 
in: The Andreas Floer Memorial Volume,
H. Hofer, C. Taubes, A. Weinstein, and
E. Zehnder, eds. 
\jour Progress in Mathematics 
\vol 133
\yr 1995
\pages 627--635
\publ Birkh\"auser
\publaddr Boston $\cdot$ Basel $\cdot$ Berlin
\endref

\ref \no \weylbook
\by H. Weyl
\book The classical groups
\publ Princeton University  Press
\publaddr Princeton NJ
\yr 1946
\endref

\ref \no \whitnone
\by H. Whitney
\paper Analytic extensions of differentiable functions defined
on closed sets
\jour Trans. Amer. Math. Soc.
\vol 36
\yr 1934
\pages  63--89
\endref

\ref \no \whitnboo
\by H. Whitney
\book Complex analytic varieties
\publ Addison-Wesley Pub. Comp.
\publaddr Reading, Ma, Menlo Park, Ca, London, Don Mills, Ontario
\yr 1972
\endref

\enddocument